\documentclass[fleqn,12pt]{article}
\usepackage{amsmath}
\usepackage{amsfonts}
\usepackage{graphicx}

\title{\Large\bf Differential geometry of line of curvature on parametric hypersurfaces in $\mathbb{E}^4$}
\author{Fatih \c{C}EL\.{I}K, Mustafa D\"ULD\"UL}
\date{\small \it Y{\i}ld{\i}z Technical University, Science and Arts Faculty, Department of Mathematics, \.Istanbul, Turkey.\\ e-mail: clk.fa17@gmail.com, mduldul@yahoo.com}

\usepackage{mathptmx}

\begin{document}
\maketitle
\begin{quote}
{\bf Abstract}
The purpose of this paper is, first, to give an algorithm that enables to obtain the lines of curvature on parametric hypersurfaces in Euclidean 4-space, and then, to obtain the curvatures of such lines by using the extended Darboux frame along the curve.
\end{quote}
{\bf 1. Introduction}\\

A curve whose velocity vectors at its every point correspond to the principal directions of the surface is called the line of curvature. Lines of curvature on surfaces have always on-going attention not only in differential geometry (see e.g. \cite{dC}, \cite{TW}) but also in geometric modeling (see e.g. \cite{Ma3}). Differential geometrical properties of such curves on parametric surfaces can be found in \cite{dC}, \cite{TW}, \cite{Ma2} and for parametric hypersurfaces in \cite{Am}. It is known that a line of curvature on a parametric surface in Euclidean 3-space $\mathbb{E}^3$ satisfies the following differential equation \cite{dC}
\begin{equation*}
(\mathsf{L}\mathsf{E}-\mathsf{N}\mathsf{F})(u')^2+(\mathsf{M}\mathsf{E}-\mathsf{N}\mathsf{G})u'v'+(\mathsf{M}\mathsf{F}-\mathsf{L}\mathsf{G})(v')^2=0.
\end{equation*}
If the above differential equation can be solved explicitly, then the line of curvature on the given surface can be obtained. In this case, it is easy to find all Frenet apparatus of the obtained line of curvature. When the equation gives the approximate solution for the line of curvature, then it is needed to construct new techniques to calculate the curvatures and Frenet vectors of the line of curvature. In 2014, Maekawa et al. present algorithms for computing the differential geometry properties of lines of curvature of parametric surfaces in $\mathbb{E}^3$. They derive the unit tangent vector, curvature vector, binormal vector, torsion, and also algorithms for evaluating the higher-order derivatives of lines of curvature of parametric surfaces \cite{Ma1}. In 2007, Che et al. studied the line of curvature and their differential geometric properties for implicit surfaces in $\mathbb{E}^3$ \cite{Ch} (The previous studies including the applications of lines of curvature have been reviewed in \cite{Ma1} and \cite{Ch}). 

Besides, lines of curvatures have also been studied in $\mathbb{E}^4$. The differential equation of the lines of curvature for immersions of surfaces into $\mathbb{E}^4$ has been established by \cite{Gu}. In \cite{So1}, the authors establish the geometric structure of the lines of curvature of a hypersurface immersed in $\mathbb{E}^4$ in a neighborhood of the set of its principal curvature singularities, consisting of the points at which at least two principal curvatures are equal. The geometric structures of the lines of curvature and the partially umbilic singularities of the three-dimensional non-compact generic quadric hypersurfaces of $\mathbb{E}^4$ have been studied by \cite{So2}.

In this paper, we give a method to compute the lines of curvature of parametric hypersurfaces in $\mathbb{E}^4$. By using the extended Darboux frame along a hypersurface curve, we also present the algorithm to obtain the Frenet vectors and the curvatures of the obtained lines of curvature.

This paper is organized as follows: In section 2, some necessary notations and reviews of the differential geometry of curves on hypersurfaces are introduced. The extended Darboux frame is also given in section 2.
In section 3, we give the method to compute the lines of curvature of parametric hypersurfaces in $\mathbb{E}^4$. The curvatures of the line of curvature are obtained in section 4.\\
\\
{\bf 2. Preliminaries}\setcounter{section}{2}\\
{\bf\footnotesize 2.1. Curves on hypersurfaces in $\mathbb{E}^4$ }\\
{\bf Definition 2.1.} The {\it ternary product} of the vectors ${\bf x}=\sum\limits_{i=1}^4 x_i{\bf e_i}$,
${\bf y}=\sum\limits_{i=1}^4 y_i{\bf e_i}$, and ${\bf z}=\sum\limits_{i=1}^4 z_i{\bf e_i}$ is
defined by \cite{SW}
\begin{equation*}
{\bf x}\otimes {\bf y}\otimes {\bf z}=\left| \begin{array}{cccc}
{\bf e}_1&{\bf e}_2&{\bf e}_3&{\bf e}_4\\
x_1      &x_2      &x_3      &x_4\\
y_1      &y_2      &y_3      &y_4\\
z_1      &z_2      &z_3      &z_4
\end{array} \right|,
\end{equation*}
where $\{{\bf e}_1, {\bf e}_2, {\bf e}_3, {\bf e}_4\}$ denotes the standard basis of $\mathbb{R}^4$.

Let  $M\subset \mathbb{E}^4$ be a regular hypersurface given by ${\bf R}={\bf R}(u_1, u_2, u_3)$ and $\alpha:I\subset \mathbb{R}\rightarrow M$ be an unit speed curve. If $\{{\bf t}, {\bf n}, {\bf b}_1, {\bf b}_2\}$ denotes the moving Frenet frame along $\alpha$, then the Frenet formulas are given by \cite{O}
\begin{equation}
{\bf t}'=k_1{\bf n},\quad\quad {\bf n}'=-k_1{\bf t}+k_2{\bf b}_1,\quad\quad {\bf b}_1'=-k_2{\bf n}+k_3{\bf b}_2,\quad\quad {\bf b}_2'=-k_3{\bf b}_1,
\end{equation}
where ${\bf t}, {\bf n}, {\bf b}_1$, and ${\bf b}_2$ denote the tangent,
the principal normal, the first binormal, and the second binormal vector fields; $k_i, (i=1,2,3)$ are the {\it i}th curvature functions of the curve $\alpha$. The Frenet vectors and curvatures of the curve are given by \cite{O}
\begin{equation}
{\bf t}=\alpha', \quad{\bf n}=\frac{\alpha''}{\|\alpha''\|}, \quad{\bf b}_2=\frac{\alpha'\otimes \alpha'' \otimes \alpha'''}{\|\alpha'\otimes \alpha'' \otimes \alpha'''\|}, \quad {\bf b}_1={\bf b}_2\otimes {\bf t}\otimes {\bf n},
\end{equation}
\begin{equation}
k_1=\|\alpha''\|,\quad k_2=\frac{\langle{\bf b}_1, \alpha'''\rangle}{k_1},\quad k_3=\frac{\langle{\bf b}_2, \alpha^{(4)}\rangle}{k_1k_2}.
\end{equation} \\

The derivatives of the curve $\alpha$ is obtained as
\begin{equation}
\alpha'={\bf t},\quad\alpha''={\bf t}^\prime=k_1{\bf n},\quad
\alpha'''=-k_1^2{\bf t}+k_1'{\bf n}+k_1k_2{\bf b}_1,
\end{equation}
\begin{equation}
\alpha^{(4)}=-3k_1k_1'{\bf t}+(-k_1^3+k_1''-k_1k_2^2){\bf n}+(2k_1' k_2
+k_1k_2'){\bf b}_1+k_1k_2k_3{\bf b}_2.
\end{equation}
In addition, we can write $\alpha'(s)=\sum\limits_{i=1}^3 {\bf R}_iu_i'$, where ${\bf R}_i=\frac{\partial {\bf R}}{\partial u_i},\,i=1,2,3$.
Also, applying the chain rule to $\alpha'$, we get the second and third derivative of the curve $\alpha$ as follows:
\begin{equation}
\alpha''(s)=\sum\limits_{i=1}^3 {\bf R}_iu_i''+\sum\limits_{i,j=1}^3 {\bf R}_{ij}u_i'u_j'
\end{equation}
\begin{equation}
\alpha'''(s)=\sum\limits_{i=1}^3 {\bf R}_iu_i'''+3\sum\limits_{i,j=1}^3 {\bf R}_{ij}u_i'u_j''+\sum\limits_{i,j,k=1}^3 {\bf R}_{ijk}u_i'u_j'u_k'.
\end{equation} 
 {\bf\footnotesize 2.2. Extended Darboux frame in $\mathbb{E}^4$ }\\

Let  $M\subset \mathbb{E}^4$ be a regular hypersurface given by ${\bf R}={\bf R}(u_1, u_2, u_3)$ and $\alpha:I\subset \mathbb{R}\rightarrow M$ be a unit speed curve. Let ${\bf T}$ denotes the unit tangent vector field along $\alpha$, and ${\bf N}$ denotes the unit normal vector field of $M$ restricted to the curve $\alpha$ given by $\displaystyle{\bf N}=\frac{{\bf R}_1\otimes {\bf R}_2\otimes {\bf R}_3}{\|{\bf R}_1\otimes {\bf R}_2\otimes {\bf R}_3\|}$. Then the extended Darboux frame of first kind along $\alpha$  is given by $\{{\bf T}, {\bf E}, {\bf D}, {\bf N}\}$, where $\displaystyle{\bf E}=\frac{{\bf T}'-\langle{\bf T}',  {\bf N}\rangle{\bf N}}{\|{\bf T}'-\langle{\bf T}',  {\bf N}\rangle{\bf N}\|}$  and ${\bf D}={\bf N}\otimes {\bf T}\otimes {\bf E}$, and it satisfies the following differential equations \cite{Du}
\begin{equation} 
\begin{pmatrix} {\bf T}'\\{\bf E}'\\{\bf D}'\\{\bf N}' \end{pmatrix}
=\begin{pmatrix} 0&k_g^{1}&0&k_n\\ -k_g^{1}&0&k_g^{2}&\tau_g^{1}\\ 0&-k_g^{2}&0&-\tau_g^{2}\\-k_n&-\tau_g^{1}&-\tau_g^{2}&0 \end{pmatrix}
 \begin{pmatrix} {\bf T} \\ {\bf E}\\ {\bf D}\\{\bf N} \end{pmatrix},
\end{equation} 
where $k_g^{i}$ and $\tau_g^{i}$ are the geodesic curvature and geodesic torsion of order $i$, respectively, $k_n$ is the normal curvature of the hypersurface in the direction of the tangent vector ${\bf T}$. Then we have \cite{Du}
 \begin{equation}
k_g^1=\langle{\bf T}',  {\bf E}\rangle, \quad k_g^2=\langle{\bf E}',  {\bf D}\rangle, \quad\tau_g^1=\langle{\bf E}',  {\bf N}\rangle, \quad\tau_g^2=\langle{\bf D}',  {\bf N}\rangle, \quad k_n=\langle{\bf T}',  {\bf N}\rangle.
\end{equation}          
\\                              
{\bf 3. Differential geometry of line of curvature of parametric hypersurfaces in $\mathbb{E}^4$}\setcounter{section}{3}\setcounter{equation}{0}\\
{\bf\footnotesize 3.1. Normal curvatures of hypersurfaces in $\mathbb{E}^4$}\\

Let ${\bf R}={\bf R}(u_1, u_2, u_3)$ denotes a regular hypersurface $M$ defined on a domain $B$. The unit normal vector field of $M$ is then $\displaystyle {\bf N}=\frac{{\bf R}_1\otimes {\bf R}_2\otimes {\bf R}_3}{\|{\bf R}_1\otimes {\bf R}_2\otimes {\bf R}_3\|}$. Since the normal curvature of $M$ at a point is given by $\displaystyle k_n=\frac{II}{I}$, we have
\begin{equation}
k_n(\lambda, \mu)=\frac{h_{11}+2h_{12}\lambda+2h_{13}\mu+2h_{23}\lambda \mu+h_{22}\lambda^2+h_{33}\mu^2}{g_{11}+2g_{12}\lambda+2g_{13}\mu+2g_{23}\lambda \mu+g_{22}\lambda^2+g_{33}\mu^2},
\end{equation}
where $\displaystyle(\lambda,\mu)=\left(\frac{du_2}{du_1},\frac{du_3}{du_1}\right)$ yields the tangent direction, $g_{ij}, h_{ij}$ denote the coefficients of the first and second fundamental forms, respectively. It is well-known that the extremal values of the normal curvature are principal curvatures, \cite{Am}. If we take the partial derivatives of $k_n$ with respect to $\lambda$ and $\mu$, we have
\begin{equation*}
\frac{\partial k_n}{\partial\lambda}=\frac{\Big(2h_{12}+2h_{23}\mu+2h_{22}\lambda\Big)I-II\Big(2g_{12}+2g_{23}\mu+2g_{22}\lambda\Big)}{\Big(g_{11}+2g_{12}\lambda+2g_{13}\mu+2g_{23}\lambda \mu+g_{22}\lambda^2+g_{33}\mu^2\Big)^2}
\end{equation*}
\begin{equation*}
\frac{\partial k_n}{\partial\mu}=\frac{\Big(2h_{13}+2h_{23}\mu+2h_{33}\lambda\Big)I-II\Big(2g_{13}+2g_{23}\mu+2g_{33}\lambda\Big)}{\Big(g_{11}+2g_{12}\lambda+2g_{13}\mu+2g_{23}\lambda \mu+g_{22}\lambda^2+g_{33}\mu^2\Big)^2}
\end{equation*}
Then,  we obtain 
\begin{equation}
k_n(\lambda, \mu)=\frac{II}{I}=\frac{h_{12}+h_{22}\lambda+h_{23}\mu}{g_{12}+g_{22}\lambda+g_{23}\mu}= \frac{h_{13}+h_{23}\lambda+h_{33}\mu}{g_{13}+g_{23}\lambda+g_{33}\mu}=\frac{h_{11}+h_{12}\lambda+h_{13}\mu}{g_{11}+g_{12}\lambda+g_{13}\mu}.
\end{equation}
Thus, the principal curvatures satisfy the following homogeneous system \cite{Am}
\begin{equation}
\left\{
\begin{array}{l}
(h_{11}-k_ng_{11})du_1+(h_{12}-k_ng_{12})du_2+(h_{13}-k_ng_{13})du_3=0\\
(h_{12}-k_ng_{12})du_1+(h_{22}-k_ng_{22})du_2+(h_{23}-k_ng_{23})du_3=0\\
(h_{13}-k_ng_{13})du_1+(h_{23}-k_ng_{23})du_2+(h_{33}-k_ng_{33})du_3=0.
\end{array}\right.
\end{equation}
Let us denote the coefficient matrix of above system with $A$, i.e.
\begin{equation*}
A=\left( \begin{array}{cccc}
 h_{11}-k_ng_{11}      & h_{12}-k_ng_{12}     &h_{13}-k_ng_{13}\\
h_{12}-k_ng_{12}      &h_{22}-k_ng_{22}      &h_{23}-k_ng_{23}\\
h_{13}-k_ng_{13}&h_{23}-k_ng_{23}&h_{33}-k_ng_{33}  
\end{array} \right).
\end{equation*}      
In the case of $rankA=0$, since all directions satisfy (3.3), the point becomes an umbilical point. Then at an umbilical point we have 
\begin{equation*}
k_n=\frac{h_{11}}{g_{11}}=\frac{h_{12}}{g_{12}}=\frac{h_{13}}{g_{13}}=\frac{h_{23}}{g_{23}}=\frac{h_{22}}{g_{22}}=\frac{h_{33}}{g_{33}}.
\end{equation*}
In the case of $rankA=3$, the system has trivial solution. For the principal directions, we need the nontrivial solutions of this system. This system has a nontrivial solution if and only if $detA=0$,
i.e.
\begin{equation*}
\begin{array}{l}
detA=\Big(g_{12}^2g_{33}+g_{13}^2g_{22}+g_{23}^2g_{11}-2g_{12}g_{23}g_{13}-g_{11}g_{22}g_{33}\Big)k_n^3\\
+\Big(2h_{12}g_{13}g_{23}-2h_{12}g_{12}g_{33}+2h_{13}g_{12}g_{23}+2h_{23}g_{12}g_{13}+h_{22}g_{11}g_{33}-h_{23}g_{11}g_{22}\\
+h_{11}g_{22}g_{33}-2h_{13}g_{22}g_{13}-2h_{23}g_{11}g_{23}-h_{23}g_{12}^2-h_{22}g_{13}^2-h_{11}g_{23}^2\Big)k_n^2\\
+\Big(2h_{12}h_{33}g_{12}-2h_{12}h_{23}g_{13}-2h_{13}h_{23}g_{12}+2h_{11}h_{23}g_{23}+2h_{13}h_{22}g_{13}-h_{11}h_{22}g_{33}-h_{22}h_{33}\\
+g_{11}-h_{11}h_{33}g_{22}-h_{12}h_{13}g_{23}+h_{23}^2g_{11}+h_{22}g_{13}^2+h_{12}^2g_{33}\Big)k_n\\
+2h_{12}h_{13}h_{23}+h_{11}h_{22}h_{33}-h_{12}^2h_{33}-h_{13}^2h_{22}-h_{11}h_{23}^2=0.
\end{array} 
\end{equation*} 
If we denote
\begin{equation*}
K_1=\frac{2h_{12}h_{13}h_{23}+h_{11}h_{22}h_{33}-h_{12}^2h_{33}-h_{13}^2h_{22}-h_{11}h_{23}^2}{g_{12}^2g_{33}+g_{13}^2g_{22}+g_{23}^2g_{11}-2g_{12}g_{23}g_{13}-g_{11}g_{22}g_{33}},
\end{equation*}
\begin{align*}
K_2=&
\frac{2h_{12}g_{13}g_{23}-2h_{12}g_{12}g_{33}+2h_{13}g_{12}g_{23}+2h_{23}g_{12}g_{13}+h_{22}g_{11}g_{33}-h_{23}g_{11}g_{22}}{3\Big(g_{12}^2g_{33}+g_{13}^2g_{22}+g_{23}^2g_{11}-2g_{12}g_{23}g_{13}-g_{11}g_{22}g_{33}\Big)}\\
&+\frac{h_{11}g_{22}g_{33}-2h_{13}g_{22}g_{13}-2h_{23}g_{11}g_{23}-h_{23}g_{12}^2-h_{22}g_{13}^2-h_{11}g_{23}^2}{3\Big(g_{12}^2g_{33}+g_{13}^2g_{22}+g_{23}^2g_{11}-2g_{12}g_{23}g_{13}-g_{11}g_{22}g_{33}\Big)},
\end{align*}

\begin{align*}
K_3=&\frac{2h_{12}h_{33}g_{12}-2h_{12}h_{23}g_{13}-2h_{13}h_{23}g_{12}+2h_{11}h_{23}g_{23}+2h_{13}h_{22}g_{13}-h_{11}h_{22}g_{33}}{3\Big(g_{12}^2g_{33}+g_{13}^2g_{22}+g_{23}^2g_{11}-2g_{12}g_{23}g_{13}-g_{11}g_{22}g_{33}\Big)}\\
&+\frac{-h_{22}h_{33}g_{11}-h_{11}h_{33}g_{22}-h_{12}h_{13}g_{23}+h_{23}^2g_{11}+h_{22}g_{13}^2+h_{12}^2g_{33}}{3\Big(g_{12}^2g_{33}+g_{13}^2g_{22}+g_{23}^2g_{11}-2g_{12}g_{23}g_{13}-g_{11}g_{22}g_{33}\Big)},
\end{align*}
we can write
\begin{equation}
k_n^3+3K_2k_n^2+3K_3k_n+K_1=0,
\end{equation}
where $K_1$ and $K_2$ correspond to the Gauss and Mean curvatures, respectively. This is a third order equation with respect to $k_n$. G. Cardano has shown the solutions for such equations \cite{EW}. By using his method, we can express the principal curvatures in terms of $K_1$, $K_2$, and $K_3$. \\
\\
{\bf\footnotesize 3.2. Computation of line of curvature }\\

Let us now give the computation technique for obtaining the line of curvature of the hypersurface. We assume that the line of curvature is arc-length parametrized. 

Let us assume that $rankA=2$. In this case, if we denote
\begin{equation}
a_1=(h_{12}-k_ng_{12})(h_{23}-k_ng_{23})-(h_{22}-k_ng_{22})(h_{13}-k_ng_{13}),
\end{equation}
\begin{equation}
a_2=(h_{12}-k_ng_{12})(h_{13}-k_ng_{13})-(h_{11}-k_ng_{11})(h_{23}-k_ng_{23}),
\end{equation}
\begin{equation}
a_3=(h_{11}-k_ng_{11})(h_{22}-k_ng_{22})-(h_{12}-k_ng_{12})^2
\end{equation}
and choose
\begin{equation}
u_1'=\eta a_1,\quad u_2'=\eta a_2,\quad u_3'=\eta a_3,
\end{equation}
(where $\eta$ is a nonzero factor) it is easy to see that (3.8) satisfies (3.3). Since the line of curvature is unit speed, its first fundamental form is given by
\begin{equation}
\sum_{i,j=1}^3 g_{ij}u_i'u_j'=1.
\end{equation}
Hence, substituting (3.8) into (3.9) determines $\eta$ as
\begin{equation}
\eta=\mp \frac{1}{\sqrt{g_{11}a_1^2+2g_{12}a_1a_2+2g_{13}a_1a_3+2g_{23}a_2a_3+g_{22}a_2^2+g_{33}a_3^2}}
\end{equation}
If we substitute (3.10) into (3.8), we obtain
\begin{equation}
\left\{\begin{array}{l}
\displaystyle
u_1'=\mp \frac{a_1}{\sqrt{g_{11}a_1^2+2g_{12}a_1a_2+2g_{13}a_1a_3+2g_{23}a_2a_3+g_{22}a_2^2+g_{33}a_3^2}},
\\
\displaystyle
u_2'=\mp \frac{a_2}{\sqrt{g_{11}a_1^2+2g_{12}a_1a_2+2g_{13}a_1a_3+2g_{23}a_2a_3+g_{22}a_2^2+g_{33}a_3^2}},
\\
\displaystyle
u_3'=\mp \frac{a_3}{\sqrt{g_{11}a_1^2+2g_{12}a_1a_2+2g_{13}a_1a_3+2g_{23}a_2a_3+g_{22}a_2^2+g_{33}a_3^2}}.
\end{array}\right.
\end{equation}
Note that the right hand sides of (3.5)-(3.7) correspond to the determinants of some submatrices obtained from the matrix $A$. For that reason, the numerators in (3.11) may be all zero. Since we assume $rankA=2$, we have a possibility in which at least one numerator for (3.11) is nonzero. %(see the Appendix for the other cases). 
Thus, the lines of curvature are then the solution of initial value problem of the triplet non-linear differential equation (3.11).
\\
\\
{\bf 4. Curvatures of line of curvature}\setcounter{section}{4}\setcounter{equation}{0}\\
{\bf\footnotesize 4.1. First curvature $(k_1)$}\\

Our aim is now to obtain all curvatures of the line of curvature $\alpha(s)$ on a parametric hypersurface. For this purpose, we need to compute the higher order derivatives of line of curvature. 

If we use (2.6) and (2.8), for the second derivative of $\alpha$ we may write
\begin{equation}
\alpha''(s)=k_g^1{\bf E}+k_n {\bf N}={\bf R}_1u_1''+{\bf R}_2u_2''+{\bf R}_3u_3''+\Omega_1,
\end{equation}
where 
$\Omega_1=\sum\limits_{i,j=1}^3 {\bf R}_{ij}u_i'u_j'$.
Since $u_1'$, $u_2'$ and $u_3'$ are known from $(3.11)$, $\Omega_1$ is known. 
By taking inner product of (3.12) with ${\bf R}_1$, ${\bf R}_2$, and ${\bf R}_3$, respectively, gives
\begin{equation}
k_g^1\langle{\bf E},  {\bf R}_1\rangle =g_{11}u_1''+g_{12}u_2''+g_{13}u_3''+\langle \Omega_1,  {\bf R}_1\rangle 
\end{equation}
\begin{equation}
k_g^1\langle{\bf E},  {\bf R}_2\rangle =g_{12}u_1''+g_{22}u_2''+g_{23}u_3''+\langle\Omega_1,  {\bf R}_2\rangle 
\end{equation}
\begin{equation}
k_g^1\langle{\bf E},  {\bf R}_3\rangle =g_{13}u_1''+g_{23}u_2''+g_{33}u_3''+\langle\Omega_1,  {\bf R}_3\rangle .
\end{equation}
Moreover, since $\alpha'=\sum\limits_{i=1}^3{\bf R}_iu_i'$ and $\langle{\bf E},  \alpha'\rangle=0$, for $u_1'\neq 0$ we have
\begin{equation}
\langle{\bf E},  {\bf R}_1\rangle =-\frac{u_2'\langle{\bf E},  {\bf R}_2\rangle }{u_1'}-\frac{u_3'\langle{\bf E},  {\bf R}_3\rangle}{u_1'}.
\end{equation}
If we substitute (4.5) into (4.2), and use (4.3), (4.4) we get
\begin{equation}
\left(\sum\limits_{i=1}^3 u_i'g_{1i}\right)u_1''+\left(\sum\limits_{i=1}^3 u_i'g_{2i}\right)u_2''+\left(\sum\limits_{i=1}^3 u_i'g_{3i}\right)u_3''=-\Big\langle\Omega_1, \alpha'\Big\rangle.
\end{equation}
By taking the derivative of the first and second equation of $(3.3)$, we obtain 
\begin{equation}
(h_{11}-k_ng_{11})u_1''+(h_{12}-k_ng_{12})u_2''+(h_{13}-k_ng_{13})u_3''=\rho_1,
\end{equation}
\begin{equation}
(h_{12}-k_ng_{12})u_1''+(h_{22}-k_ng_{22})u_2''+(h_{23}-k_ng_{23})u_3''=\rho_2,
\end{equation}
where
\begin{align*}
 \rho_1&=-\Big(h_{11}'-k_n'g_{11}-k_ng_{11}'\Big)u_1'-\Big(h_{12}'-k_n'g_{12}-k_ng_{12}'\Big)u_2'-\Big(h_{13}'-k_n'g_{13}-k_ng_{13}'\Big)u_3',\\
 \rho_2&=-\Big(h_{12}'-k_n'g_{12}-k_ng_{12}'\Big)u_1'-\Big(h_{22}'-k_n'g_{22}-k_ng_{22}'\Big)u_2'-\Big(h_{23}'-k_n'g_{23}-k_ng_{23}'\Big)u_3'
\end{align*}
 By combining (4.3), (4.4), (4.6), (4.7) and (4.8), we obtain the following nonhomogeneous linear equation system
\begin{equation} 
\begin{pmatrix}
\sum\limits_{i=1}^3 u_i'g_{1i}&\sum\limits_{i=1}^3 u_i'g_{2i}&\sum\limits_{i=1}^3 u_i'g_{3i}&0&0\\
g_{12}&g_{22}&g_{23}&-1&0\\g_{13}&g_{23}&g_{33}&0&-1\\ h_{11}-k_ng_{11} &h_{12}-k_ng_{12}& h_{13}-k_ng_{13}&0&0 \\ h_{12}-k_ng_{12} & h_{22}-k_ng_{22} & h_{23}-k_ng_{23}&0&0 
\end{pmatrix}
\begin{pmatrix} u_1''\\u_2''\\u_3''\\k_g^1\langle{\bf E},  {\bf R}_2\rangle\\ k_g^1\langle{\bf E},  {\bf R}_3\rangle \end{pmatrix} 
= \begin{pmatrix} 
-\langle\Omega_1, \alpha'\rangle\\
-\langle\Omega_1,  {\bf R}_2\rangle\\ -\langle\Omega_1,  {\bf R}_3\rangle\\  \rho_1 \\ \rho_2 \end{pmatrix} .
\end{equation}  
Since the determinant of the coefficient matrix is
\begin{equation*}  \sqrt{g_{11}a_1^2+2g_{12}a_1a_2+2g_{13}a_1a_3+2g_{23}a_2a_3+g_{22}a_2^2+g_{33}a_3^2}\not=0,
\end{equation*} 
$u_1''$, $u_2''$, $u_3''$, $k_g^1\langle{\bf E},  {\bf R}_2\rangle$, and $k_g^1\langle{\bf E},  {\bf R}_3\rangle$  can be evaluated from $(4.9)$ which enable us to compute the curvature vector $\alpha''$, ${\bf E}=\frac{{\bf T}'-\langle{\bf T}',  {\bf N}\rangle{\bf N}}{\|{\bf T}'-\langle{\bf T}',  {\bf N}\rangle{\bf N}\|}$, and $k_g^1$. Hence, the first curvature of a line of curvature is obtained by
\begin{equation}
k_1=\sqrt{(k_n)^2+(k_g^1)^2}.
\end{equation}
{\bf\footnotesize 4.2. Second curvature $(k_2)$}\\

 To obtain the second curvature, we need to determine the third order derivative of $\alpha$.

Since $\alpha$ is a line of curvature, we have $\tau_g^1=\tau_g^2=0$. If we take the derivative of $\alpha''=k_g^1{\bf E}+k_n{\bf N}$ with respect to arc-length, we get 
\begin{equation}
\alpha'''=(k_g^1)'{\bf E}-k_1^2{\bf T}+k_g^1k_g^2{\bf D}+k_n'{\bf N}.
\end{equation}
Thus, if we use (2.8)
\begin{equation}
(k_g^1)'{\bf E}-k_1^2{\bf T}+k_g^1k_g^2{\bf D}+k_n'{\bf N}={\bf R}_1u_1'''+{\bf R}_2u_2'''+{\bf R}_3u_3'''+\Omega_2,
\end{equation}
where $\Omega_2=3\sum\limits_{i,j=1}^3 {\bf R}_{ij}u_i'u_j''+\sum\limits_{i,j,k=1}^3 {\bf R}_{ijk}u_i'u_j'u_k'$.
If we take the inner product of both sides of (4.12) with ${\bf R}_1$, ${\bf R}_2$, ${\bf R}_3$, respectively, we obtain the linear equation system
\begin{equation}
g_{11}u_1'''+g_{12}u_2'''+g_{13}u_3'''-(k_g^1)'\langle{\bf E},  {\bf R}_1\rangle-k_g^1k_g^2\langle{\bf D},  {\bf R}_1\rangle=-k_1^2\langle{\bf T},  {\bf R}_1\rangle-\langle\Omega_2,  {\bf R}_1\rangle
\end{equation}
\begin{equation}
g_{12}u_1'''+g_{22}u_2'''+g_{23}u_3'''-(k_g^1)'\langle{\bf E},  {\bf R}_2\rangle-k_g^1k_g^2\langle{\bf D},  {\bf R}_2\rangle=-k_1^2\langle{\bf T},  {\bf R}_2\rangle-\langle\Omega_2,  {\bf R}_2\rangle
\end{equation}
\begin{equation}
g_{13}u_1'''+g_{23}u_2'''+g_{33}u_3'''-(k_g^1)'\langle{\bf E},  {\bf R}_3\rangle-k_g^1k_g^2\langle{\bf D},  {\bf R}_3\rangle=-k_1^2\langle{\bf T},  {\bf R}_3\rangle-\langle\Omega_2,  {\bf R}_3\rangle.
\end{equation}
in which $u_1'''$, $u_2'''$, $u_3'''$, $(k_g^1)'$ and $k_g^2$ are unknowns. So, we have to find two more equations to obtain these unknowns.
By taking the derivatives of $(4.7)$ and $(4.8)$, we get 
\begin{equation}
(h_{11}-k_ng_{11})u_1'''+(h_{12}-k_ng_{12})u_2'''+(h_{13}-k_ng_{13})u_3'''=\rho_3
\end{equation}
\begin{equation}
(h_{12}-k_ng_{12})u_1'''+(h_{22}-k_ng_{22})u_2'''+(h_{23}-k_ng_{23})u_3'''=\rho_4,
\end{equation}
where
\begin{equation*}
\begin{array}{ll}
 \rho_3=&-\Big(h_{11}''-k_n''g_{11}-2k_n'g_{11}'-k_ng_{11}''\Big)u_1'-2\Big(h_{11}'-k_n'g_{11}-k_ng_{11}\Big)u_1''\\
&-\Big(h_{12}''-k_n''g_{12}-2k_n'g_{12}'-k_ng_{12}''\Big)u_2'-2\Big(h_{12}'-k_n'g_{12}-k_ng_{12}\Big)u_2''\\
&-\Big(h_{13}''-k_n''g_{13}-2k_n'g_{13}'-k_ng_{13}''\Big)u_3'-2\Big(h_{13}'-k_n'g_{13}-k_ng_{13}\Big)u_3'',
\end{array}
\end{equation*}
\begin{equation*}
\begin{array}{ll}
 \rho_4=&-\Big(h_{12}''-k_n''g_{12}-2k_n'g_{12}'-k_ng_{12}''\Big)u_1'-2\Big(h_{12}'-k_n'g_{12}-k_ng_{12}\Big)u_1''\\
&-\Big(h_{22}''-k_n''g_{22}-2k_n'g_{22}'-k_ng_{22}''\Big)u_2' -2\Big(h_{22}'-k_n'g_{22}-k_ng_{22}\Big)u_2''\\
&-\Big(h_{23}''-k_n''g_{23}-2k_n'g_{23}'-k_ng_{23}''\Big)u_3'-2\Big(h_{23}'-k_n'g_{23}-k_ng_{23}\Big)u_3''.
\end{array}
\end{equation*}

By combining (4.13)-(4.17), we get the nonhomogeneous linear equation system $QX=S$, where
\begin{equation*} 
Q=\begin{pmatrix}g_{11}&g_{12}&g_{13}&-\langle{\bf E},  {\bf R}_1\rangle&-k_g^1\langle{\bf D},  {\bf R}_1\rangle\\g_{12}&g_{22}&g_{23}&-\langle{\bf E},  {\bf R}_2\rangle&-k_g^1\langle{\bf D},  {\bf R}_2\rangle\\ g_{13}&g_{23}&g_{33}&-\langle{\bf E},  {\bf R}_3\rangle&-k_g^1\langle{\bf D},  {\bf R}_3\rangle\\h_{11}-k_ng_{11} &h_{12}-k_ng_{12}& h_{13}-k_ng_{13}&0&0 \\ h_{12}-k_ng_{12} & h_{22}-k_ng_{22} & h_{23}-k_ng_{23}&0&0 \end{pmatrix}
\end{equation*}
\begin{equation*} 
X=
\begin{pmatrix}
u_1'''\\u_2'''\\u_3'''\\(k_g^1)'\\ k_g^2 
\end{pmatrix}, 
\quad S= \begin{pmatrix} -k_1^2\langle{\bf T},  {\bf R}_2\rangle-\langle\Omega_2,  {\bf R}_2\rangle \\-k_1^2\langle{\bf T},  {\bf R}_2\rangle-\langle\Omega_2,  {\bf R}_2\rangle\\-k_1^2\langle{\bf T},  {\bf R}_3\rangle-\langle\Omega_2,  {\bf R}_3\rangle\\ \rho_3\\ \rho_4 
\end{pmatrix}. 
\end{equation*}
The unknowns $u_1'''$, $u_2'''$, $u_3'''$, $(k_g^1)'$  and $k_g^2$ can be obtained when $detQ\not=0$. Thus, we find $\alpha'''$ which enables us to compute the Frenet vectors ${\bf b}_1, {\bf b}_2$ by using (2.2). Therefore, the second curvature can be obtained by $\displaystyle k_2= \frac{\langle{\bf b}_1,  \alpha'''\rangle}{k_1}$. Also, $ k_1'=\langle{\bf n}, \alpha'''\rangle. $\\
\\
{\bf\footnotesize 4.3. Third curvature $(k_3)$}\\

Similarly, for the fourth derivative of the line of curvature we may write
\begin{align}
\alpha^{(4)}&=-3k_1k_1'{\bf T}+(-k_1^3+k_1''-k_1k_2^2){\bf n}+(2k_1'k_2+k_1k_2'){\bf b}_1+ k_1k_2k_3{\bf b}_2 \nonumber \\
&={\bf R}_1u_1^{(4)}+{\bf R}_2u_2^{(4)}+{\bf R}_3u_3^{(4)}+\Omega_3,
\end{align}
where 
\begin{equation*}
\Omega_3=4\sum\limits_{i,j=1}^3 {\bf R}_{ij}u_i'''u_j'+3\sum\limits_{i,j=1}^3 {\bf R}_{ij}u_i''u_j''+6\sum\limits_{i,j,k=1}^3 {\bf R}_{ijk}u_i'u_j'u_k'+\sum\limits_{i,j,k,\ell=1}^3 {\bf R}_{ijk\ell}u_i'u_j'u_k'u_\ell'.
\end{equation*}
If we take the inner product of the both sides of (4.18) with ${\bf R}_1$, ${\bf R}_2$ and ${\bf R}_3$, respectively, we get the following equations: 
\begin{equation}
\begin{array}{l}
{g_{11}u_1^{(4)}+g_{12}u_2^{(4)}+g_{13}u_3^{(4)}-k_1''\langle{\bf n},  {\bf R}_1\rangle -k_1k_2'\langle{\bf b}_1,  {\bf R}_1\rangle-k_1k_2k_3\langle{\bf b}_2,  {\bf R}_1\rangle} \\ 
=-\langle \Omega_3, {\bf R}_1\rangle-3k_1k_1'\langle{\bf T},  {\bf R}_1\rangle-(k_1^3+k_1k_2^2)\langle{\bf n},  {\bf R}_1\rangle+2k_1'k_2\langle{\bf b}_1,  {\bf R}_1\rangle,
\end{array}
\end{equation}
\begin{eqnarray}
\lefteqn{g_{12}u_1^{(4)}+g_{22}u_2^{(4)}+g_{23}u_3^{(4)}-k_1''\langle{\bf n},  {\bf R}_2\rangle-k_1k_2'\langle{\bf b}_1,  {\bf R}_2\rangle-k_1k_2k_3\langle{\bf b}_2,  {\bf R}_2\rangle} \nonumber \\
&=-\langle \Omega_3, {\bf R}_2\rangle-3k_1k_1'\langle{\bf T},  {\bf R}_2\rangle-(k_1^3+k_1k_2^2)\langle{\bf n},  {\bf R}_2\rangle+2k_1'k_2\langle{\bf b}_1,  {\bf R}_2\rangle,
\end{eqnarray}
\begin{eqnarray}
\lefteqn{g_{13}u_1^{(4)}+g_{23}u_2^{(4)}+g_{33}u_3^{(4)}-k_1''\langle{\bf n},  {\bf R}_3\rangle-k_1k_2'\langle{\bf b}_1,  {\bf R}_3\rangle-k_1k_2k_3\langle{\bf b}_2,  {\bf R}_3\rangle} \nonumber \\
&=-\langle \Omega_3, {\bf R}_3\rangle-3k_1k_1'\langle{\bf T},  {\bf R}_3\rangle-(k_1^3+k_1k_2^2)\langle{\bf n},  {\bf R}_3\rangle+2k_1'k_2\langle{\bf b}_1,  {\bf R}_3\rangle.
\end{eqnarray}
If we take the inner product of $(4.18)$ with ${\bf n}$, ${\bf b}_1$ and ${\bf b}_2$, respectively, we obtain
\begin{equation}
\langle{\bf n}, {\bf R}_1\rangle u_1^{(4)}+\langle{\bf n}, {\bf R}_2\rangle u_2^{(4)}+\langle{\bf n}, {\bf R}_3\rangle u_3^{(4)}-k_1''=-\langle \Omega_3, {\bf n}\rangle-k_1^3-k_1k_2^2
\end{equation}
\begin{equation}
\langle{\bf b}_1, {\bf R}_1\rangle u_1^{(4)}+\langle{\bf b}_1, {\bf R}_2\rangle u_2^{(4)}+\langle{\bf b}_1, {\bf R}_3\rangle u_3^{(4)}-k_1k_2'=-\langle \Omega_3, {\bf b}_1\rangle+2k_1'k_2
\end{equation}
\begin{equation}
\langle{\bf b}_2, {\bf R}_1\rangle u_1^{(4)}+\langle{\bf b}_2, {\bf R}_2\rangle u_2^{(4)}+\langle{\bf b}_2, {\bf R}_3\rangle u_3^{(4)}-k_1k_2k_3=-\langle \Omega_3, {\bf b}_2\rangle.
\end{equation}
Thus, (4.19)-(4.24) constitute a linear equation system with the unknowns $u_1^{(4)}$, $u_2^{(4)}$, $u_3^{(4)}$, $k_1''$, $k_2'$, and the third curvature $k_3$. These unknowns can be computed when $detM\not=0$, where

\begin{eqnarray} 
M=\begin{pmatrix}
g_{11}&g_{12}&g_{13}&-\langle{\bf n},  {\bf R}_1\rangle&-k_1\langle{\bf b}_1,  {\bf R}_1\rangle&-k_1k_2\langle{\bf b}_2,  {\bf R}_1\rangle\\g_{12}&g_{22}&g_{23}&-\langle{\bf n},  {\bf R}_2\rangle&-k_1\langle{\bf b}_1,  {\bf R}_2\rangle&-k_1k_2\langle{\bf b}_2,  {\bf R}_2\rangle\\ g_{13}&g_{23}&g_{33}&-\langle{\bf n},  {\bf R}_3\rangle&-k_1\langle{\bf b}_1,  {\bf R}_3\rangle&-k_1k_2\langle{\bf b}_2,  {\bf R}_3\rangle\\ \langle{\bf n},  {\bf R}_1\rangle &\langle{\bf n},  {\bf R}_2\rangle&\langle{\bf n},  {\bf R}_3\rangle&-1&0&0 \\ \langle{\bf b}_1,  {\bf R}_1\rangle & \langle{\bf b}_1,  {\bf R}_2\rangle & \langle{\bf b}_1,  {\bf R}_3\rangle&0&-k_1&0\\\langle{\bf b}_2,  {\bf R}_1\rangle&\langle{\bf b}_2,  {\bf R}_2\rangle &\langle{\bf b}_2,  {\bf R}_3\rangle&0&0&-k_1k_2
\end{pmatrix}.
\end{eqnarray}


\begin{thebibliography}{99}

\bibitem{Ma1}  H.K. Joo, T. Yazaki, M. Takezawa, T. Maekawa,  Differential geometry properties of lines of curvature of parametric surfaces and their visualition, \emph{Graphical Methods}, \textbf{76} (2014), 224--238.

\bibitem{Ma2}  N.M. Patrikalakis, T. Maekawa, \emph{Shape Interrogation for Computer Aided Design Manufacturing}, Springer-Verlag, Heidelberg, 2002. 

\bibitem{Ma3}  T. Maekawa, F.-E. Wolter, N.M. Patrikalakis, Umbilics and lines of curvature for shape interrogation, \emph{Comp. Aided. Geom. Design}, {\bf 13} (1996), 133--161. 

\bibitem{Ch} W. Che, J. Paul, Z.X. Zhang, Lines of curvature and umbilic points for implicit surfaces, \emph{Comp. Aided. Geom. Design} {\bf 24} (2007), 395--409.

\bibitem{dC} M. P. do Carmo, \emph{Differential Geometry of Curves and Surfaces}, Prentice Hall, Englewood Cliffs, NJ, 1976. 

\bibitem{Du} M. D{\"{u}}ld{\"{u}}l, B. Uyar D{\"{u}}ld{\"{u}}l, N. Kuruo\u{g}lu, E. {\"{O}}zdamar, Extension of the Darboux frame into Euclidean $4$-space and its invariants, \emph{Turk. J. Math.} \textbf{41} (2017),  1628--1639.

\bibitem{O}
O. Alessio, Differential geometry of intersection curves in $\mathbb{R}^4$ of three implicit surfaces, \emph{Comp. Aided. Geom. Design} \textbf{26} (2009), 455--475. 

\bibitem{Am} Yu. Aminov, 
 \emph{The Geometry of Submanifolds}, Gordon and Breach Science Publishers, 2001.

\bibitem{TW} T.J. Willmore, \emph{An Introduction to Differential Geometry}, Oxford University Press, 1959.

\bibitem{SW} M. Z. Williams, F. M. Stein, A triple product of vectors in four space, 
\emph{Mathematics Magazine}, \textbf{37} (1964), 230--235.

\bibitem{Gu}  C. Gutierrez, I. Guadalupe, R. Tribuzy, V. Gu\'i\~nez, A Differential Equation for Lines of Curvature on Surfaces Immersed in $\mathbb{R}^4$, \emph{Qual. Theory Dyn. Syst.}, {\bf 2} (2001), 207--220. 

\bibitem{So1}  D. Lopes, J. Sotomayor, R. Garcia, Partially umbilic singularities of hypersurfaces of $\mathbb{R}^4$, \emph{Bull. Sci. math.}, {\bf 139} (2015), 431--472.

 \bibitem{So2} J. Sotomayor, R. Garcia, Lines of Curvature on Quadric Hypersurfaces of $\mathbb{R}^4$, \emph{Lobachevskii Journal of Mathematics}, {\bf 37}(3) (2016), 288--306.

\bibitem{EW} E.W. Weisstein, Cubic formula, http://mathworld.wolfram.com/CubicFormula.html

\end{thebibliography}
\end{document}